%% file: main_eng.tex
\documentclass[11pt, oneside, a4paper]{amsart}
\usepackage{amsfonts, amscd, amsmath, amssymb}
\usepackage{hyperref}
\usepackage[english]{babel}

\author{Sergey Maksymenko}
\title[Sections of Lie group actions]
{Sections of Lie group actions \\ and a theorem by M.~Newman}
\curraddr{Topology dept., Institute of Mathematics of NAS of Ukraine, 
Te\-re\-shchen\-kivs\-ka str. 3, Kiyv, 01601, Ukraine}
\email{maks@imath.kiev.ua}
\keywords{Lie group, fixed point, flow}
\subjclass[2000] {
37C10, 
37C27, 
57S15, 
}

\input{macros.tex}

\begin{document}

\begin{abstract}
Let $\manif$ be a smooth finite-dimensional manifold, $\grp$ a Lie group, and $\gact:\grp \times \manif\to\manif$ a smooth action.
Consider the following mapping
$$ \actsect: \smsgrp \to \smmap, $$
defined by $\actsect(\afunc)(x) = \actnote{\afunc(x)}{x}$,
for $\afunc\in\smsgrp$ and $x\in\manif$.
In this paper we describe the structure of inverse images of elements of $\smmap$
under $\actsect$ for $\dim\grp=1$, i.e. when $\grp$ is either $\RRR$ or $S^1$.
As an application we obtain a new proof of the well-known theorem by M.~Newman concerning the interior of the fixed point set of a Lie group action.

This paper is a translation from russian of the paper published in
Proceedings of the conference ``Fundamental mathematics today'' devoted to the 10-th annivesary of the Independent University of Moscow, 2003, p.~246--258.
\end{abstract}

 \maketitle
 \input{01intro.tex}
 \input{02prop.tex}
 \input{03regp.tex}

 \input{04perlm.tex}
 \input{05fixp.tex}
 \input{06prfr1.tex}
 \input{07prfs1.tex}
 \input{08bibl.tex}

\end{document}

%% file: macros.tex
\def\globcnt{subsection}

\swapnumbers
\theoremstyle{plain}

\newtheorem{theorem}[\globcnt]{Theorem}
\newtheorem{lemma}[\globcnt]{Lemma}
\newtheorem{prop}[\globcnt]{Proposition}
\newtheorem{cor}[\globcnt]{Corollary}
\theoremstyle{definition}
\newtheorem{defn}[\globcnt]{Definition}

\newtheorem{examp}[\globcnt]{Example}

\makeatletter
\@addtoreset{subsubsection}{section}
\@addtoreset{subsubsection}{subsection}
\makeatother

\newcommand\ZidRTheorem{th:Zid_descr_R1}
\newcommand\ZidSTheorem{th:Zid_descr_S1}
\newcommand\ZidCompTheorem{th:compG_IntFix_empty}

\theoremstyle{plain}

\providecommand{\eqref}[1]{(\ref{#1})}
\renewcommand{\eqref}[1]{(\ref{#1})}

\newcommand\CCC{{\mathbb C}}
\newcommand\NNN{{\mathbb N}}

\newcommand\RRR{{\mathbb R}}

\newcommand\ZZZ{{\mathbb Z}}

\newcommand\FF{{\mathcal F}}

\newcommand\JJ{{\mathcal J}}

\newcommand\OO{{\mathcal O}}

\newcommand\Int{\mathrm{Int}\,}

\newcommand\id{\mathrm{id}}

\newcommand\Fr{\mathrm{Fr}\,}

\newcommand\IM{\mathrm{im}\,}

\newcommand\diag{\mathrm{diag}\,}

\newcommand\Fix{{\mathrm{Fix}}\,}
\newcommand\Per{{\mathrm{Per}}\,}

\newcommand\manif{M}

\def\Diff{\mathrm{Diff}{}}             

\def\Fld{\FF}

\def\flow{\Phi}

\newcommand\smsp[1]{C^{\infty}\ifx#1\empty\else(#1)\fi}
\newcommand\smfunc{\smsp{\manif,\RRR}}
\newcommand\smmap{\smsp{\manif,\manif}}

\newcommand\afunc{\alpha}
\newcommand\tafunc{\tilde\alpha}

\def\amap{f}

\newcommand\bfunc{\beta}
\newcommand\bmap{g}

\newcommand\cfunc{\gamma}
\newcommand\cmap{h}

\newcommand\ufunc{\sigma}

\newcommand\grp{G}            
\newcommand\gact{\flow}      
\newcommand\kergact{K}      
\newcommand\tkergact{\widetilde{\kergact}}
\newcommand\actsect{\varphi}  

\newcommand\pnt{z}

\newcommand\Orbpnt{\OO_{\pnt}}

\newcommand\smsgrp{\smsp{\manif,\grp}}

\newcommand\actnote[2]{{#1}\cdot{#2}}
\newcommand\actnotethree[3]{\actnote{\actnote{#1}{#2}}{#3}}

\newcommand\Zid{Z_{\id}}

\newcommand\nbb{V}
\newcommand\pntb{y}
\newcommand\nbc{W}
\newcommand\pntc{z}

\newcommand\sectflow{\smsp{\manif, \RRR}}

\newcommand\rcomp{C}
\newcommand\nbhcomp{V}
\newcommand\nind{\lambda}
\newcommand\GL{GL}

\newcommand\Rot[1]{R\left(#1\right)}
\newcommand\dflow{\Psi}
\newcommand\expflow{\nu}

\newcommand\dimM{m}

\newcommand\tgact{\widetilde{\gact}}
\newcommand\tZid{\widetilde{\Zid}}

\newcommand\Odzflow{U}
\newcommand\Interval{\JJ}
\newcommand\Shift{\actsect}

\newcommand\sectlocfl{\smsp{\Odzflow, \Interval}}
\newcommand\shiftlocfl{\smsp{\Odzflow,\manif}}

\newcommand\tlambda{\widetilde{\lambda}}
\newcommand\NPerSet{N}
\newcommand\PerSet{P}
\newcommand\RegSet{R}

\newcommand\Zzero{\Zid'}
\newcommand\Fixgact{F}
\newcommand\Fixflow{\Fixgact}
\newcommand\manifInt{\manif\setminus\Int\Fixflow}
\newcommand\evalhom{\tau_{\pnt}}
\newcommand\boldz{\mathbf{z}}

\newcommand\cmp{c}
\newcommand\rl{d}

%% file: 01intro.tex
\section{Introduction}
Let $\manif$ be a smooth ($C^{\infty}$) connected finite-dimensional manifold, $\grp$ a Lie group and 
\begin{equation}\label{equ:gact}
\gact:\grp \times \manif\to\manif
\end{equation}
a smooth action of $\grp$ on $\manif$.
We will often write $\actnote{\bmap}{\pnt}$ instead $\gact(\bmap,\pnt)$
for $\bmap\in\grp$ and $\pnt\in\manif$.

For each point $\pnt\in\manif$ let $\Orbpnt \subset \manif$ be the corresponding orbit of $\pnt$.
Let also $\Fix\gact$ designates the fixed-point set of $\gact$.
A point which is not fixed will be called \emph{regular}.

Suppose that a mapping $\amap:\manif\to\manif$ preserves each orbit of $\gact$, i.e. $\amap(\Orbpnt) \subset \Orbpnt$ for all $\pnt\in\manif$.
The following question arises often in different problems: can $\amap$ be ``smoothly parametrized'' by elements of $\grp$?
Thus we need a smooth mapping $\afunc:\manif\to\grp$ such that
$$
 \amap(\pnt)=\gact(\afunc(\pnt), \pnt) = \actnote{\afunc(\pnt)}{\pnt} ,
$$
for all $\pnt\in\manif$?
This question leads to the consideration of the following mapping:
\begin{equation}\label{equ:actsect}
\actsect: \smsgrp \to \smmap,
\end{equation}
defined by 
$$
\actsect(\afunc)(\pnt) = \gact(\afunc(\pnt),\pnt) = \actnote{\afunc(\pnt)}{\pnt},
$$
for $\afunc\in\smsgrp$ and $\pnt\in\manif$.

Evidently, the identity mapping $\id_{\manif}$ of $\manif$ belongs to the image of $\actsect$.
Consider the inverse image $\Zid = \gact^{-1}(\id_{\manif}) \subset \smsgrp$.
This set of smooth mappings $\manif\to\grp$ plays an important role for the understanding the structure of $\actsect$.
It is easy to see that $\Zid$ is a subgroup of the $\smsgrp$.
Moreover, we will show that the inverse images of mappings $\smmap$ of $\actsect$
coincide with the adjacent classes of $\smsgrp$ by $\Zid$.

The main result of this paper is a full description of the group $\Zid$ for the case 
$\dim\grp = 1$, i.e. when $\grp$ is either $\RRR^1$ or $S^1$.

For simplicity let us denote by $\Fixgact$ the fixed-point set 
$\Fix\gact$ of the action~\eqref{equ:gact}.
Let also $\Int\Fixgact$ and $\Fr(\Int\Fixgact)$ be respectively the interior of $\Fixgact$ in $\manif$ and the frame of this interior.

\begin{theorem}\label\ZidRTheorem
Suppose that $\grp=\RRR$, i.e. $\gact$ is a flow on $\manif$.

\emph{(1)} If $\Int\Fixgact \not= \emptyset$, then
$$\Zid = \{ \mu \in \sectflow \ | \ \mu|_{\manif \setminus \Int\Fixgact} = 0\}.$$

\emph{(2)} Suppose that $\Int\Fixgact = \emptyset$. 
Then we have two possibility: either $\Zid=\{0\}$ or
there exists a smooth strictly positive smooth function $\mu:\manif\to(0,\infty)$,
such that
$$\Zid = \{ n\mu \ | \ n\in\ZZZ \} \approx \ZZZ.$$
In this case each regular point $\pnt$ of $\flow$ is periodic and $\mu(\pnt)$ is equal to the period of $\pnt$.
\end{theorem}

\begin{theorem}\label\ZidSTheorem
Suppose $\grp=S^1$. 
Let $\kergact$ be the ineffectivity of the action~\eqref{equ:gact},
i.e. the kernel of the induced homomorphism $\grp \to \Diff\manif$.
If the action $\gact$ is non-trivial, then $\Int\Fixgact = \emptyset$ and $\Zid$ consists of constant mappings $\manif\to\kergact \subset \grp$.
Consequently, $\Zid$ is isomorphic to $\kergact$ and is a finite cyclic group.
\end{theorem}

As an application we obtain a new proof of the following theorem~\ref\ZidCompTheorem.
It is a variant of  the well-known theorem by M.~Newman~\cite{Newman} concerning the actions of cyclic groups, see also~\cite{Montgomery57, MontSamZip} for more general account.

\begin{theorem}\label\ZidCompTheorem
Let $\grp$ be a compact Lie group acting non-trivially on a finite-dimensional manifold.
Then the set of fixed points of $\grp$ is nowhere dense.
\end{theorem}
\proof
Since $\grp$ is compact, it contains at least one one-parametric subgroup isomorphic to $SO(2)\approx S^1$.
Indeed, let $\grp_1$ be arbitrary $1$-parametric subgroup in $\grp$.
Then its closure $\overline{\grp_1}$ is a compact abelian Lie group, i.e. a torus, and therefore has $S^1$-subgroups.

Let $\Fix S$ be the fixed-point set of the induced action of $S$ on $\manif$.
Then $\Fix\grp \subset \Fix S$.
By Theorem~\ref\ZidSTheorem, $\Fix S$ is nowhere dense, therefore so is $\Fix\grp$.
\endproof

\subsection{Structure of the paper.}
In \S\ref{sect:GenProp} we consider the general properties of the mapping $\actsect$ that hold for arbitrary Lie groups $\grp$.
Further, we shall confine ourselves with the case $\grp=\RRR$ only.
In \S\ref{sect:RegPnt} we consider the behavior of functions belonging to $\Zid$ near regular points of $\flow$. 
\S \ref{sect:LinMapPer} includes a lemma about the lower bound of the periods of trajectories of linear flows.
This lemma will be used in \S\ref{sect:FixPerPoints} where we prove two statements about local behavior of functions from $\Zid$.
Finally, in \S\ref{sect:ZidRProof} and \S\ref{sect:ZidSProof} we prove Theorems~\ref\ZidRTheorem and~\ref\ZidSTheorem.

%% file: 02prop.tex
\section{Properties of $\actsect$.}\label{sect:GenProp}

\begin{lemma}
The image of $\IM\actsect$ is a subsemigroup in $\smmap$.
Moreover, the intersection $\IM\actsect \cap \Diff\manif$ is a subgroup of $\Diff\manif$.
\end{lemma}
\proof
Let $\afunc, \bfunc, \cfunc \in \smsgrp$ and let $\amap,\bmap, \cmap \in \smmap$ be respectively their images under $\actsect$.
Suppose also that $\cmap$ is a diffeomorphism.
For the proof of lemma we have to find mappings 
$$
\ufunc_{\amap\circ\bmap},
\ufunc_{\cmap^{-1}}:\manif\to\grp,
$$
such that 
$\amap\circ\bmap = \actsect(\ufunc_{\amap\circ\bmap})$
and
$\cmap^{-1} = \actsect(\ufunc_{\cmap^{-1}})$.
It is not easy to see that the following mappings satisfy the conditions above:
\begin{gather}
\label{equ:f_g_shift}
\ufunc_{\amap\circ\bmap}(\pnt) =
\actnote{\afunc(\bmap(\pnt))} {\bfunc(\pnt)}, \\
\label{equ:inv_h_shift}
\ufunc_{\cmap^{-1}}(\pnt) =
\left( \cfunc \left( \cmap^{-1}(\pnt) \right)\right)^{-1},
\end{gather}
for all $\pnt\in\manif$.
Indeed,
$$
\amap \circ \bmap(\pnt) =
\actnote{\afunc(\bmap(\pnt))} {\bmap(\pnt)} =
\actnotethree{ \afunc(\bmap(\pnt))}{\bfunc(\pnt)} {\pnt} =
\actnote{\ufunc_{\amap\circ\bmap}(\pnt)}{\pnt}.
$$
For the proof of~\eqref{equ:inv_h_shift} notice that the identity
$\cmap \left( \cmap^{-1}(\pnt) \right)=\pnt$ means that 
$ \actnote{ \cfunc \left(\cmap^{-1}(\pnt) \right)}{\cmap^{-1}(\pnt)} = \pnt$,
whence,
$$
\cmap^{-1}(\pnt) =
\actnote{ \cfunc \left( \cmap^{-1}(\pnt) \right)^{-1}}{\pnt} =
\actnote{\ufunc_{\cmap^{-1}}(\pnt)}{\pnt}.\qed
$$

\subsection{}
Denote by $\Zid(\gact)$ the inverse image of the identity mapping:
$\actsect^{-1}(\id_{\manif}) \subset \smsgrp$
\begin{equation}\label{equ:Zid_for_acitons}
\Zid(\gact) := \actsect^{-1}(\id_{\manif}).
\end{equation}
Then $\actnote{\mu(\pnt)}{\pnt} = \pnt$ for all $\mu\in\Zid(\gact)$ и $\pnt\in\manif$.

The following statement is evident.
\begin{prop}
The set $\Zid = \Zid(\gact)$ has the following properties:
\begin{enumerate}
\item
\label{enum:Zid_subgr}
$\Zid$ is a subgroup in$\smsgrp$.

\item
\label{enum:ab_1_in_Zid}
Let $\afunc,\bfunc\in \smsgrp$.
Then $\actsect(\afunc)=\actsect(\bfunc)$ if and only if 
$\actnote{\afunc^{-1}}{\bfunc} \in \Zid$. 
\qed
\end{enumerate}
\end{prop}

\begin{lemma}\label{lm:x_in_IntFix}
Suppose that $\Int\Fixgact\not=\emptyset$ and let $\afunc,\bfunc\in\smsgrp$ coincide outside the interior of $\Int\Fixgact$:
$$\afunc|_{\manif\setminus\Int\Fixgact} = \bfunc|_{\manif\setminus\Int\Fixgact}.$$
Then
$\actsect(\afunc)=\actsect(\bfunc)$.
In particular, if $\afunc(\pnt) = e \in\grp$ for all $\pnt\in \manif\setminus\Int\Fixgact$, then $\afunc\in\Zid$.
\end{lemma}
\proof
We have to prove that $\actnote{\afunc(\pnt)}{\pnt} = \actnote{\bfunc(\pnt)}{\pnt}$ for $\pnt\in\manif$.
Let $\pnt\in\manif\setminus\Int\Fixgact$.
Then $\afunc(\pnt)=\bfunc(\pnt)$, whence $\actnote{\afunc(\pnt)}{\pnt} = \actnote{\bfunc(\pnt)}{\pnt}.$

Suppose that $\pnt\in\Int\Fixgact$.
Then $\actnote{t}{\pnt} = \pnt$ for each $t \in \grp$, whence $\actnote{\afunc(\pnt)}{\pnt} = \actnote{\bfunc(\pnt)}{\pnt} = \pnt.$
\endproof

The following lemma is also evident.
\begin{lemma}\label{lm:mu_const_g_in_ker}
A constant mapping $\mu:\manif\to\grp$ belongs to $\Zid$ if and only if the image of $\mu$ belongs to the ineffectivity kernel of the action $\gact$.
\qed
\end{lemma}

%% file: 03regp.tex
\section{Regular points of flows}\label{sect:RegPnt}
In the sequel, it is assumed that $\grp=\RRR$.
It is also convenient to assume that $\gact$ is a \emph{local\/} action.
Thus $\gact$ is a local flow defined on some open subset of $\manif$.
Let us recall the definitions.

\begin{defn}
Let $\Odzflow$ be an open connected subset of $\manif$ and $\Interval$ an open interval in $\RRR$ containing $0$.
A smooth mapping
\begin{equation}\label{equ:local_flow}
\flow:\Interval\times\Odzflow \to \manif
\end{equation}
is a {\em local flow} if the following conditions hold true:
\begin{enumerate}
\item[(1)] $\flow(0,x) = x, \forall x\in\Odzflow$,
\item[(2)] $\flow(s, \flow(t,x)) = \flow(t+s,x)$,
for $x\in\Odzflow$ and $t,s\in\Interval$ provided $\flow(t,x)\in\Odzflow$ and $t+s\in\Interval$.
\end{enumerate}

In the case $\Odzflow=\manif$ and $\Interval=\RRR$ the flow $\flow$ is {\em global}.
For each $t\in\Interval$ the restriction of $\flow$ to $\{t\} \times\Odzflow$
$$\flow|_{\{t\} \times\Odzflow}:\Odzflow \to \manif$$
will be denoted by $\flow_t$.

Let $\pnt\in\manif$.
The {\em orbit\/} of $\pnt$ is the following set $\flow(\Interval \times \{x\}) \subset \manif$.
A point $\pnt\in\Odzflow$ is {\em fixed} with respect to the flow provided
$\flow(t,\pnt)=\pnt$ for each $t\in\Interval$.
All other points are {\em regular}.
A regular point $\pnt$ is {\em periodical} if $\flow(t,\pnt)=\pnt$ for some $t>0$.
The least such number $t$ is called the {\em period\/} of $\pnt$ and is denoted by $\Per(\pnt)$.
The orbit of a periodic point is {\em closed}, and the orbit of a regular but non-periodic is {\em non-closed}.
\end{defn}

\subsection{}
Notice that a local flow~\eqref{equ:local_flow} induces the following mapping
\begin{equation}\label{equ:loc_Shift}
\Shift:\sectlocfl \to \shiftlocfl
\end{equation}
defined by $\Shift(\afunc)(\pnt) = \flow(\afunc(\pnt),\pnt)$, for $\pnt\in\Odzflow$ and $\afunc\in\sectlocfl$.

Let $\id_{\Odzflow}:\Odzflow \to \manif$ be the identity embedding and designate
$$\Zid := \Shift^{-1}(\id_{\Odzflow}).$$

\begin{lemma}\label{lm:mu_const_traj}
Let $\omega$ be a non-constant orbit of $\flow$ and $\mu\in\Zid$.
If $\omega$ is non-closed, then $\mu|_{\omega}=0$.
If $\omega$ is closed orbit of period $\theta$, then $\mu|_{\omega}= n \theta$ for some $n \in \ZZZ$.
\end{lemma}
\proof
Recall that the condition $\mu\in\Zid$ means that $\flow(\mu(x),x) = x$ for all $x\in\manif$.

Consider a non-closed orbit $\omega$ of $\flow$.
Then for arbitrary pair $x,y\in\omega$ there exists a {\em unique\/} number $t \in \Interval$ such that $\flow(t,x) = y$.
In particular, $t=0$ iff $x=y$.
Therefore $\mu(x)= 0$ for all $x\in\omega$.

Let $\omega$ be a non-closed orbit of period $\theta$ and $x\in\omega$.
Then $\flow(t,x)=x$ if and only if $t = n \theta$ for some $n\in\ZZZ$.
Therefore $\mu(x)=n(x)\theta$ for all $x\in\omega$,
where $n = \mu / \theta:\omega\to\ZZZ$ is a {\em continuous} function.
It follows that $n$ is constant, i.e. $\mu|_{\omega}=n\theta$.
\endproof

The following lemma claims a local uniqueness of function of $\Zid$
in a neighborhood of a regular point of a flow.

\begin{lemma}\label{lm:reg_per_comp}
Let $\rcomp$ be a connected component of the set of regular points of $\flow$.
Let also $\afunc,\bfunc\in\sectlocfl$ be such that $\Shift(\afunc)=\Shift(\bfunc)$.
If $\afunc(\pntb)=\bfunc(\pntb)$ for some $\pntb\in\rcomp$,
then $\afunc|_{\rcomp} = \bfunc|_{\rcomp}$.
In particular, if $\afunc\in\Zid$ and $\afunc(\pntb)=0$, then $\afunc|_{\rcomp}=0$.
\end{lemma}
\proof
It suffices to show that $\afunc=\bfunc$ in a neighborhood of $\pntb\in\rcomp$.
We will now express a function $\afunc$ through $\Shift(\afunc)$ in a neighborhood of a regular point of $\flow$.

Let $\amap=\Shift(\afunc)$ and $a=\afunc(\pntb)$.
Since the point $\pntc = \amap(\pntb)=\flow_{a}(\pntb)$ is also regular, there are local coordinates $(x_1,\ldots,x_n)$ in a some neighborhood $\nbc$ of $\pntc$ such that $\pntc=0$ and
$\flow(t, x_1,\ldots,x_n) = (x_1 + t, x_2,\ldots,x_n)$.
Consider the following neighborhood $\nbb = \amap^{-1}(\nbc) \cap \flow_{-a}(\nbc)$ of $\pntb$. Evidently,
\begin{equation}\label{equ:reg_pnt_lsect}
  \afunc(x) = p_1\circ\Shift(\afunc)\circ\flow(a,x) - p_1\circ\flow(a,x), \ \ \forall x \in\nbb.
\end{equation}
where $p_1:\RRR^n \to \RRR$ is a projection onto the first coordinate.

This formula proves our lemma.
Indeed, if $\Shift(\afunc)=\Shift(\bfunc)$ and $\afunc(\pntb)=\bfunc(\pntb)=a$,
then by formula~\eqref{equ:reg_pnt_lsect} $\afunc\equiv\bfunc$ in some neighborhood of $\pntb$ consisting of regular points only.
Therefore they coincides on $\rcomp$.
\endproof

%% file: 04perlm.tex
\section{Periods of linear flows}\label{sect:LinMapPer}
We prove here a lemma that gives a lower bounds for the periods of orbits of linear flows.
First we introduce some notation and recall the ``real'' Jordan form of a matrix.

Let $A$ be a $(k \times k)$-matrix.
Then a {\em Jordan cell \/} $J_p(A)$ of $A$ is the following $(p k \times p k)$-matrix:
$$
J_p(A) =
\left.
\left\|
\begin{matrix}
A     & 0      & \cdots & 0 & 0 \\
E_k   & A      & \cdots & 0 & 0\\
\cdot & \cdot  & \cdot  & \cdot & \cdot \\
0 & 0 & \cdots & A      & 0 \\
0 & 0 & \cdots & E_k    & A
\end{matrix}
\right\|
\ \right\} p
,
$$
where $E_k$ is a unit $(k \times k)$-matrix.
For $\alpha, \beta\in\RRR$ set
\begin{equation}\label{equ:rotation_matrix}
\Rot{\alpha + i\beta} = \left( \begin{matrix}\alpha & -\beta \\ \beta & \alpha \end{matrix}\right).
\end{equation}

Let $X$ be a square matrix with real entries and $\lambda\in\CCC\setminus\RRR$ an eigen value of $X$. 
Then $\overline{\lambda}$ is also an eigen value of $X$.
Let $\lambda_1, \overline{\lambda_1}, \ldots, \lambda_{\cmp}, \overline{\lambda_{\cmp}}$ be all complex and $\lambda_{\cmp+1},\ldots,\lambda_{\cmp+\rl}$ all real eigen values of $X$.
Then by Jordan theorem $X$ is conjugated to the following matrix:
\begin{equation}\label{equ:Real_JordForm}
\diag[
J_{p_1}(\Rot{\lambda_{1}}), \ldots, J_{p_s}(\Rot{\lambda_{\cmp}}),
J_{k_1}(\lambda_{\cmp+1}), \ldots, J_{k_m}(\lambda_{\cmp+\rl})].
\end{equation}
(see e.g. Theorem~2.2.5 in~\cite{PalisMelo}.)

We will also need the following formula for the exponent of a Jordan cell:
\begin{equation}\label{equ:exp_J_A}
e^{J_{k}(A) t} =
\left\|
\begin{array}{cccc}
e^{At}                               & 0                                    & \cdots  & 0      \\
t \cdot e^{At}                       & e^{At}                               & \cdots  & 0      \\
\cdots                               & \cdots                               & \cdots  & \cdots \\
\frac{t^{k-1}}{(k-1)!}  \cdot e^{At} & \frac{t^{k-2}}{(k-2)!}  \cdot e^{At} & \cdots  & e^{At}
\end{array}
\right\|.
\end{equation}

Let $M(n)$ be the algebra of all real $n \times n$-matrices and $\exp:M(n)\to\GL(\RRR,n)$ an exponential mapping.
\begin{lemma}\label{lm:linear_map_periods}
Let $A \in M(n)$ and let $\Lambda = \{ \lambda_j \}_{j=1}^{r}$ be the set of its purely imagine non-zero eigen values.
Then the linear flow $\flow(t,x)=e^{At}x$ has a closed orbit if and only if $\Lambda \not= \emptyset$.
In this case the period of each closed orbit of $\flow$ is $\geq\min\limits_{j=1..r} \frac{2\pi}{|\lambda_j|} $.
\end{lemma}
\proof
First suppose that $A$ is a Jordan cell.
Then $\flow$ has closed orbits iff $A = J^p(\Rot{i\beta})$ for some $\beta\in\RRR\setminus\{0\}$.
In this case all eigen values of $A$ are equal to $\pm i\beta$.
Moreover, it follows from formula~\eqref{equ:exp_J_A} that all closed orbits of $\flow$ belong to the invariant subspace generated by two latter coordinates.
They also have the same period $\frac{2\pi}{|\beta|}$ = $\frac{2\pi}{|\lambda|}$.

Consider now the general case.
We can assume that $A$ has a real normal Jordan form~\eqref{equ:Real_JordForm},
and that for $i=1..r$ $(r\leq\cmp)$ the eigen values $\lambda_i$ constitute $\Lambda$.
Set $m=\cmp+\rl$ and designate by $V_i \subset \RRR^n (i=1..m)$ the invariant subspace of $\flow $ corresponding to the corresponding cell either $J_{p_i}(\lambda_i)$ or $J_{p_i}(\Rot{\lambda_{i}})$.
Then $\RRR^n = \mathop\oplus\limits_{i=1}^{m} V_i$.

Let $p_i:\RRR^n \to V_i$ be a natural projection and $\flow^i = \flow|_{V_i}$ the restriction of $\flow$ to $V_i$.
Then for every orbit $\omega$ of $\flow$ and every $i=1..m$ the set $\omega_i = p_i(\omega)$
is the orbit of the flow $\flow^i$.
Moreover, $\omega$ is closed iff all $\omega_i$ are either closed of constant and for at least one $j=1..r$ the $\omega_j$ is closed.
Therefore, $\flow$ has a closed orbit iff $\Lambda\not=\emptyset$.

Suppose that $\omega$ is a closed orbit of $\flow$ of a period $\theta$.
Let $\omega_j = p_j(\omega)$ be the projection of $\omega$ being a closed orbit of some $\flow^j$.
Then the period of $\omega_j$ is equal to $\theta_j=\frac{2\pi}{|\lambda_j|}$.
Since the projection $p_j$ factors $\flow$ onto $\flow^j$, i.e. $p_j\circ\flow_t=\flow^j_t\circ p_j$, it follows that $\theta = s \theta_j$ for some $s\in\NNN$.
In particular, $\theta \geq \theta_j \geq \min\limits_{k=1..r} \frac{2\pi}{|\lambda_k|}$.
\endproof

\begin{cor}\label{cor:linear_map_periods}
Let $\{ A_i \}_{i\in\NNN} \subset M(n)$ be a sequence of matrices such that for all $i\in\NNN$ the linear flow $\flow_i(t,x) = e^{A_i t}x$ has closed orbits.
Let $\theta_i$ be the minimum of periods of orbits of $\flow_i$.
If $\lim\limits_{i \to \infty} A_i = 0$, then $\lim\limits_{i \to \infty} \theta_i = \infty$.
\end{cor}
\proof
Let $\Lambda_{i}$ be the set of those eigen values of $A_i$ that correspond to closed orbits of $\flow_i$ (see Lemma~\ref{lm:linear_map_periods}.)
Then $\Lambda_{i}\not=\emptyset$ for all $i\in\NNN$.
Set $\tlambda_i = \max\limits_{\lambda \in\Lambda_i} |\lambda|$.
Then $\theta_i = 2\pi / \tlambda_i$.
Since $\lim\limits_{i \to \infty} A_i = 0$, it follows from the continuity of a spectrum of matrices that $\lim\limits_{i \to \infty} \tlambda_i = 0$.
Therefore by Lemma~\ref{lm:linear_map_periods}, $\lim\limits_{i\to\infty}
\theta_i = \lim\limits_{i\to\infty} 2\pi / \tlambda_i =\infty$.
\endproof

%% file: 05fixp.tex
\section{Fixed and periodic poitns of flows}\label{sect:FixPerPoints}
\begin{prop}\label{pr:x_in_FrFix}
Let $\RegSet = \Odzflow\setminus\Fixflow$ be the set of regular points of a flow $\flow$,
$\nbhcomp$ be a connected component of $\RegSet$, and $\pnt\in\Fixflow\cap\overline{\nbhcomp}$.
Let also $\mu\in\Zid$ be such that $\mu(\pnt)=0$.
Then $\mu \equiv 0$ on $\overline{\nbhcomp}$.
\end{prop}
\proof
Notice that the components of $\RegSet$ can be divided into the following two parts
$$ \RegSet = \NPerSet \cup \PerSet, $$
where $\NPerSet$ consists of those components that include at least one non-closed orbit of $\flow$ and $\PerSet$ consists of components with only periodic orbits. 

Let $\mu\in\Zid$.
Then it follows from Lemmas~\ref{lm:mu_const_traj} and~\ref{lm:reg_per_comp} that 
$\mu|_{\overline{\NPerSet}} = 0$.
Thus it suffices to prove our proposition for the component belonging to $\nbhcomp \subset \PerSet$.

Let $\{ \pnt_i\}_{i\in\NNN} \subset \nbhcomp$ be a sequence of periodic points of $\flow$ converging to $\pnt$.
For each $i\in\NNN$ let $\theta_i$ be the period of $\pnt_i$.
Then by Lemma~\ref{lm:mu_const_traj}, $\mu(\pnt_i) = n_i\theta_i$ for some $n_i\in\ZZZ$.
Hence by continuity of $\mu$ we get
\begin{equation}\label{equ:mu_ni_thi}
 \mu(\pnt_i) = n_i \theta_i \to \mu(\pnt) = 0.
\end{equation}

Taking if necessary a subsequence we can assume that there exists a certain finite or infinite limit $\theta = \lim\limits_{i \to \infty} \theta_i \geq 0$.
By Corollary~\ref{cor:linear_map_periods} we have that $\theta>0$.
Then it follows from~\eqref{equ:mu_ni_thi} that $n_i=0$ for all sufficiently large $i\in\NNN$.
In particular, $\mu(\pnt_i)=0$ for some $i\in\NNN$.
Since $\pnt_i\in\nbhcomp_\nind$, we get from Lemma~\ref{lm:reg_per_comp} that $\mu\equiv 0$ on $\nbhcomp_\nind$.
Proposition is proved.
\endproof

\begin{prop}\label{pr:x_in_FrIntFix}
Let $\Fr(\Fixflow) = \Fixflow\setminus\Int\Fixflow$ be the boundary of the fixed-point set of $\flow$.
Suppose that for a point $\pnt\in\Fr(\Fixflow)$ one of the following conditions holds true:
\begin{enumerate}
\item 
$\pnt$ belongs to the boundary of the interior of $F$, i.e. $\pnt\in\Fr(\Int\Fixflow)$;

\item
the tangent linear flow at $\pnt$ is trivial, i.e.
$\frac{\partial\flow}{\partial x}(t, \pnt)=E_n$ for all $t\in\Interval$.
\end{enumerate}

Then for each $\mu\in\Zid$ we have that $\mu \equiv 0$ in some neighborhood of $\pnt$ in the set $\Odzflow\setminus\Int\Fixflow$.
\end{prop}
\proof
Since the problem is local, we can assume that $\manif=\RRR^{\dimM}$ and $\pnt=0\in\RRR^{\dimM}$ is the origin.

Define the following mapping $\dflow: \Interval \times \Odzflow \to \GL(\RRR,n)$ by
$\dflow(t,x) = \frac{\partial \flow}{\partial x}(t,x)$.
Since  $\dflow(0,x)=E_n$ for all $x \in \Odzflow$, we see that the mapping 
$\expflow = \exp^{-1} \circ \dflow:\Interval \times \Odzflow \to M(n)$
is defined in some neighborhood of $(0,\pnt)$ in $\Interval \times \Odzflow $.
Thus $\dflow(t,x) = e^{\expflow(t,x)}$.

Notice that for each $x\in\Odzflow$ the restriction $\dflow(*,x):\Interval \to \GL(\RRR,n)$
is a local {\em homomorphism}. 
Therefore it induces a linear flow on $\RRR^{\dimM}$.
Hence the matrix $A(x,t) = \expflow(t,x)/t$ does not depend on $t\in\Interval$,
i.e. $\dflow(t,x) = e^{A(x) t}$.

Moreover, for each periodic point $x$ the flow $\dflow(*,x)$ has closed trajectories.
Indeed, let $\Fld(x) = \frac{\partial\flow}{\partial t}(0,x)$ be a vector field generating $\flow$.
Applying to both parts of the following relation 
$$ \flow(s,\flow(t,x)) = \flow(t,\flow(s,x)) $$
the operator $\frac{\partial}{\partial t}$ and then set $s=0$ we will obtain
$$
\begin{CD}
\frac{\partial\flow}{\partial t} (0, \flow(t,x)) =
\frac{\partial\flow}{\partial x} (t, x) \frac{\partial\flow}{\partial t} (0,x),
\end{CD}
$$
whence $\Fld(\flow(t,x)) = \dflow(t,x) \Fld(x)$.
This means that the vectors $\Fld(\flow(t,x))$ and $\Fld(x)$ belong to the same orbit of the flow $\dflow(*,x)$.
It follows that if $x$ is a periodic point of $\flow$, then $\Fld(x)$ is a periodic point of $\dflow(*,x)$ of the same period: $\Per(x) \geq \Per(\Fld(x))$.

Now we can complete the proposition.
Evidently, (2) holds for each internal point of $\Fixflow$.
Therefore it also holds for the boundary points of $\Fr(\Int\Fixflow)$.
Hence (1) implies (2).

Thus suppose that (2) holds true.
Consider the following sequences of vectors
$$F_i(t) = \frac{\partial\flow}{\partial t} (t, \pnt_i)$$
and matrices
$$A_i(t) = \frac{\partial\flow}{\partial x} (t, \pnt_i)$$
depending on a parameter $t\in\Interval$.

As noted above, since $\pnt_i$ is a periodic point for $\flow$, it follows that each vector $F_i(t)$ is also periodic of the same period $\leq \Per(\pnt_i)=\theta_i$.
By (2) $\lim\limits_{i\to\infty} A_i(t) = E_n$.
Then from Corollary~\ref{cor:linear_map_periods} we get 
$$
\theta = \lim\limits_{i\to\infty} \theta_i \geq
  \lim\limits_{i\to\infty} \Per(F_i(t)) = \infty.
$$
Since the value $\mu(\pnt) = \lim\limits_{i\to\infty} n_i\theta_i$ is finite, we see that 
$\lim\limits_{i\to\infty} n_i = 0$.
Hence $\mu(\pnt) = 0$.
\endproof

%% file: 06prfr1.tex
\section{Proof of Theorem~\ref\ZidRTheorem} \label{sect:ZidRProof}

\subsection*{Case (1)}
Suppose that $\Int\Fixflow\not=\emptyset$.
Denote
$$\Zzero := \{ \mu \in \sectflow \ | \ \mu|_{\manif \setminus \Int\Fixflow} = 0\}.$$
We should show that $\Zzero = \Zid$.
By Lemma~\ref{lm:x_in_IntFix} $\Zzero \subset \Zid$.

Let $\mu\in\Zid$.
By statement (1) of Proposition~\ref{pr:x_in_FrIntFix} we have that $\mu(\pnt)=0$ for each $\pnt\in\Fr(\Int\Fixflow)$.
Then by Proposition~\ref{pr:x_in_FrFix} $\mu=0$ on each connected component $\manif\setminus\Int\Fixflow$ containing $\pnt$.

Since $\manif$ is connected, it follows that each connected component of the set 
$\manif\setminus\Int\Fixgact$ intersects $\Fr(\Int\Fixgact)$.
Therefore $\mu=0$ on $\manif\setminus\Int\Fixgact$, i.e. $\mu\in\Zzero$.
Hence $\Zid\subset\Zzero$.
This proves (1).

\subsection*{Case (2)} 
Let $\Int\Fixgact=\emptyset$.
Suppose that $\Zid\not=\{0\}$.
We should prove that there exists a smooth function $\mu > 0$ such that $\Zid = \{ n\cdot \mu\}_{n\in\ZZZ}$.

Consider at first arbitrary $\mu\in\Zid$.
If $\mu(\pnt)=0$ for at least one point $\pnt\in\manif$, then by Proposition~\ref{pr:x_in_FrFix}  $\mu \equiv 0$ on $\manifInt=\manif$.
Thus if $\mu(\pnt)\not=0$, then $\mu\not=0$ on $\manif$.
Therefore we can assume that $\mu>0$ on all of $\manif$.

For each point $\pnt\in\manif$ define the following mapping
$\evalhom:\Zid\to\RRR$ by $\evalhom(\nu) = \nu(\pnt)$.
It is easy to see that $\evalhom$ is a {\em homomorphism.}
By Proposition~\ref{pr:x_in_FrFix} its kernel is trivial $\ker\evalhom=0$, i.e. $\evalhom$ is injective.

Let $\pnt$ be a regular point of $\flow$.
Then by Lemma~\ref{lm:mu_const_traj} we see that $\IM\evalhom$ is a closed subgroup in $\RRR$.
Hence $\IM\evalhom$ is either trivial or isomorphic with $\ZZZ$.

Suppose that $\IM\evalhom \approx \ZZZ$.
Let $r\in\IM\evalhom \subset \RRR$ be a positive generator of $\IM\evalhom$.
Then the function $\mu = \evalhom^{-1}(r)$ is a strictly positive generator of the group $\Zid$.
This proves (2) and Theorem~\ref\ZidRTheorem.
\endproof

\examp
Consider the following two flows on the complex plane:
$$
\Phi(t,\boldz) = e^{2\pi i (1+|\boldz|^2) \cdot t } \boldz,
 \qquad \text{and}\qquad
\Psi(t,\boldz) = e^{2\pi i |\boldz|^2  \cdot t} \boldz.
$$
Let us calculate the groups $\Zid(\Phi)$ and $\Zid(\Psi)$.

Evidently, the flows $\Phi$ and $\Psi$ have the same trajectories: concentric circles with the center at the origin $0 \in\CCC$.
Nevertheless the corresponding tangent linear flows at $0$ of $\Phi$ and $\Psi$ differ each from other:
$$
\frac{ \partial\Phi }{\partial\boldz} (t,0)\xi= e^{2\pi i t }\xi,
 \qquad \qquad
\frac{ \partial\Psi }{\partial\boldz} (t,0)\xi = \xi,
$$
where $\xi$ is a tangent vector at $0$.
The latter equality means that 
$\frac{ \partial\Psi }{\partial\boldz} (t,0) =
\left(\begin{smallmatrix}1&0\\0&1 \end{smallmatrix}\right)$.
Then by Propositions~\ref{pr:x_in_FrFix} and~\ref{pr:x_in_FrIntFix} $\Zid(\Psi)=\{0\}$.

On the other hand, it is easy to see that the following function $\mu(\boldz) = \frac{1}{1 + |\boldz|^2}$ belongs to $\Zid(\Phi)$.
Whence $\Zid(\Phi)\not=\{0\}$.
Since the fixed-point set $\Fix\Phi=\{0\}$ is nowhere dense in $\CCC$, it follows from Theorem~\ref\ZidRTheorem that $\Zid\approx \ZZZ$.
It is also clear that for each point $\boldz\not=0$ the value $\mu(\boldz)$ is equal to the period of $\boldz$.
Then again by Theorem~\ref\ZidRTheorem we see that $\mu$ is a generator of the group
$$\Zid(\Phi) = \{n\cdot \mu(\boldz) \}_{n\in\ZZZ}.$$
\endexamp

%% file: 07prfs1.tex
\section{Proof of Theorem~\ref\ZidSTheorem}\label{sect:ZidSProof}
Let $p:\RRR\to S^1$ be the universal covering of $S^1$.
Define the following mapping $$\tgact:\RRR\times\manif\to\manif$$
by $\tgact(t,\pnt) = \gact(p(t),\pnt)$.
Then $\tgact$ is the action of the group $\RRR$ on $\manif$, i.e. a flow, that covers the action $\gact$.
Set
$$\tZid = \Zid(\tgact) = \{ \tafunc\in\smfunc \ | \
      \actnote{\tafunc(\pnt)}{\pnt} = \pnt, \forall \pnt\in\manif \}.$$

Suppose that $\afunc:\manif\to S^1$ belongs to $\Zid$, i.e. 
$\gact(\afunc(x),x)=x$ for all $x\in\manif$.
Let $\pnt\in\manif$ be a point and $\Odzflow$ a small neighborhood of $\pnt$.
Then $\afunc$ lifts to a function $\tafunc:\Odzflow\to\RRR$ such that
$p(\tafunc(x))=\afunc(x)$ for all $x\in\Odzflow$.
Hence
$$\tgact(\tafunc(x),x)=\gact(p(\tafunc(x)),x) = \gact(\afunc(x),x) = x.$$
This means that $\tafunc$ belongs to $\tZid$ ``locally''.
Notice also that we can always choose $\tafunc$ so that $\tafunc(\pnt)\not=0$.

Suppose that $\pnt\in\Fr(\Int\Fix\gact)$.
Then by Theorem~\ref\ZidRTheorem applied to the restriction of $\flow$ to $\Odzflow$ we get that $\tafunc(\pnt)=0$.
This contradicts to the choice of $\tafunc$.
Hence such a point $\pnt$ does not exists, i.e. $\Fr(\Int\Fixgact)=\emptyset$.
This is possible only when either $\Int\Fixgact=\manif$ or $\Int\Fixgact=\emptyset$.
In the former case the action $\tgact$ is trivial, whence $\Int\Fixgact=\emptyset$.

Since $p(\ZZZ)=1 \in S^1$, we see that the group $\tZid$ includes all constant functions $\manif \to \ZZZ$.
Therefore is includes more than one element.
Whence by (2) of Theorem~\ref\ZidRTheorem \ $\tZid \approx \ZZZ$.
Therefore $\tZid$ consists of constant mappings from $\manif$ to some subgroup $P$ of $\RRR$ such that $P\approx\ZZZ$.
It follows from Lemma~\ref{lm:mu_const_g_in_ker} that $P$ coincides with the innefectivity kernel $\tkergact$ of the action $\tgact$.

Therefore the subgroup $p(\tkergact) \subset S^1$ is the ineffectivity kernel $\kergact$ of the action $\gact$ and $\Zid$ consists of constant mappings $\manif\to \kergact$.
Moreover, since the kernel of the homomorphism $p$ is infinite, we obtain that $\Zid \approx \kergact \approx \tkergact/\ker p$ is a finite cyclic group.
\endproof

%% file: 08bibl.tex
\newcommand\Bibitem[5]{
\bibitem[#1]{#2}
{\sc #3}
\emph{#4}
#5
}